\documentclass{amsart}
\include {mak}
\newtheorem{theorem}{Theorem}[section]
\newtheorem{lemma}[theorem]{Lemma}
\newtheorem{proposition}[theorem]{Proposition}
\newtheorem{corollary}[theorem]{Corollary}
\theoremstyle{definition}
\newtheorem{definition}[theorem]{Definition}

\theoremstyle{remark}

\numberwithin{equation}{section}

\begin{document}

% Your \newcommans below (if there are any):

\oddsidemargin 16.5mm
\evensidemargin 16.5mm

\thispagestyle{plain}

\begin{center}

{\Large\bf SUPERSTABILITY OF ADJOINTABLE MAPPINGS ON HILBERT
$C^*$-MODULES
\rule{0mm}{6mm}\renewcommand{\thefootnote}{}%Enter at least one, but not more than 4 MSCs.
% First entered MSC will be a primary one, others (at most 3) will be secondary.
\footnotetext{\scriptsize 2000 Mathematics Subject Classification: Primary 46L08; Secondary 47B48, 39B52, 46L05\\
Keywords and Prases: Hyers--Ulam--Rassias stability, superstability,
Hilbert $C^*$-module, $C^*$-algebra, $\varphi$-perturbation of an
adjointable mapping. }}

\vspace{1cc} {\large\it Michael Frank, Pasc G\u avru\c ta and
Mohammad Sal Moslehian}

\vspace{1cc}

{\scriptsize Dedicated to the Memory of Professor D.
S.~Mitrinovi\'c}

\vspace{1cc}

\parbox{24cc}{{\scriptsize\bf
We define the notion of $\varphi$-perturbation of a densely defined
adjointable mapping and prove that any such mapping $f$ between
Hilbert ${\mathcal A}$-modules over a fixed $C^*$-algebra ${\mathcal
A}$ with densely defined corresponding  mapping $g$ is ${\mathcal
A}$-linear and adjointable in the classical sense with adjoint $g$.
If both $f$ and $g$ are everywhere defined then they are bounded.
Our work concerns with the concept of Hyers--Ulam--Rassias stability
originated from the Th.M.~Rassias' stability theorem that appeared
in his paper [On the stability of the linear mapping in Banach
spaces, Proc. Amer. Math. Soc. 72 (1978), 297--300]. We also
indicate interesting complementary results in the case where the
Hilbert $C^*$-modules admit non-adjointable $C^*$-linear mappings.
}}
\end{center}

\vspace{1cc}

%Body of the text here. Sections should be numbered consecutively, their titles centered.
%
%References in the text should read like e.g. Baker[1]  rather than
%Baker(1971).
%
%Numbered formulae should be recalled in the text (prepared by \LaTeX) using their labels,
%rather than numbers. Formulae which are not recalled in the paper
%should not be numbered.

\vspace{1.5cc}

\section{INTRODUCTION}

We say a functional equation $(\mathcal E)$ is \textit{stable} if
any function $g$ approximately satisfying the equation $(\mathcal
E)$ is near to an exact solution of $(\mathcal E)$. The equation
$(\mathcal E)$ is called \textit{superstable} if every approximate
solution of $(\mathcal E)$ is indeed a solution (see \cite{BAK} for
another notion of superstability namely {\it superstability modulo
the bounded functions}). More than a half century ago, S.M.~Ulam
\cite{ULA} proposed the first stability problem which was partially
solved by D.H.~Hyers \cite{HYE} in the framework of Banach spaces.
Later, T.~Aoki \cite{AOK} proved the stability of the additive
mapping and Th.M.~Rassias \cite{RAS1} proved the stability of the
linear mapping for mappings $f$ from a normed space into a Banach
space such that the norm of the Cauchy difference $f(x+y) - f(x) -
f(y)$ is bounded by the expression $\varepsilon(\|x\|\sp p+\|y\|\sp
p)$ for some $\varepsilon\ge 0$, for some $ 0 \leq p < 1 $ and for
all $x, y$. The terminology ``Hyers--Ulam--Rassias stability'' was
indeed originated from Th.M.~Rassias's paper \cite{RAS1}. In 1994, a
further generalization was obtained by P. G\u avru\c ta \cite{GAV},
in which he replaced the bound $\varepsilon(\|x\|\sp p+\|y\|\sp p)$
by a general control function $\varphi(x, y)$. This terminology can
be applied to functional equations and mappings on various
generalized notions of Hilbert spaces; see \cite{AMY, A-M, B-C-M}.
We refer the interested reader to monographs \cite{CZE1, CZE2,
H-I-R, JUN, MIT, RAS3} and references therein for more information.

The notion of Hilbert $C^*$-module is a generalization of the notion
of Hilbert space. This object was first used by I.~Kaplansky
\cite{KAP}. Interacting with the theory of operator algebras and
including ideas from non-commutative geometry it progresses and
produces results and new problems attracting attention, see
\cite{FRA,LAN,M-T}.

\noindent Let ${\mathcal A}$ be a $C^*$-algebra and ${\mathcal X}$
be a complex linear space, which is a right ${\mathcal A}$-module
with a scalar multiplication satisfying $\lambda(xa)=x(\lambda
a)=(\lambda x)a$ for $x \in {\mathcal X},a \in {\mathcal A}, \lambda
\in {\mathbb C}$. The space ${\mathcal X}$ is called a (right)
pre-Hilbert ${\mathcal A}$-module if there exists an ${\mathcal
A}$-inner product $\langle .,.\rangle :{\mathcal X} \times {\mathcal
X}\to {\mathcal A}$ satisfying

(i) $\langle x,x\rangle\geq 0$ and $\langle x,x\rangle=0$~~~ if and
only if~~~ $x=0$;

(ii) $\langle x,y+\lambda z\rangle=\langle x,y\rangle+\lambda
\langle x,z\rangle$;

(iii) $\langle x,ya\rangle=\langle x,y\rangle a$;

(iv) $\langle x,y\rangle^*=\langle y,x\rangle$;

\noindent for all $x, y, z \in {\mathcal X},\, \lambda \in {\mathbb
C},\, a \in {\mathcal A}$. The pre-Hilbert module ${\mathcal X}$ is
called a (right) Hilbert ${\mathcal A}$-module if it is complete
with respect to the norm $\| x \| =\| \langle x,x\rangle\| ^{1/2}$.
Left Hilbert ${\mathcal A}$-modules can be defined in a similar way.
Two typical examples are

(I) Every inner product space is a left pre-Hilbert ${\mathbb
C}$-module.

(II) Let ${\mathcal A}$ be a $C^*$-algebra. Then every norm-closed
right ideal $I$ of ${\mathcal A}$ is a Hilbert ${\mathcal A}$-module
if one defines $\langle a, b\rangle = a^*b \quad (a, b \in I)$.

\noindent A mapping $f:{\mathcal X}\to {\mathcal Y}$ between Hilbert
${\mathcal A}$-modules is called adjointable if there exists a
mapping $g:{\mathcal Y}\to {\mathcal X}$ such that $\langle
f(x),y\rangle=\langle x,g(y)\rangle$ for all $x\in {\mathcal D}(f)
\subseteq {\mathcal X}, y\in {\mathcal D} \subseteq {\mathcal Y}$.
Throughout the paper, we assume that $f$ and $g$ are both everywhere
defined or both densely defined. The unique mapping $g$ is denoted
by $f^*$ and is called the adjoint of $f$.

An ${\mathcal A}$-linear bounded operator $K$ on a Hilbert
${\mathcal A}$-module ${\mathcal X}$ is called ``compact'' if it
belongs to the norm-closed linear span of the set of all elementary
operators $\theta_{x,y}\,\,(x,y \in {\mathcal X})$ defined by
$\theta_{x,y}(z)= x \langle y,z\rangle \,\, (z\in {\mathcal X})$.

In this paper, we prove the superstability of adjointable mappings
on Hilbert $C^*$-modules in the spirit of Hyers--Ulam--Rassias and
indicate interesting complementary results in the case where the
Hilbert $C^*$-modules admit non-adjointable $C^*$-linear mappings.

\section{MAIN RESULTS}

Throughout this section, ${\mathcal A}$ denotes a $C^*$-algebra,
${\mathcal X}$ and ${\mathcal Y}$ denote Hilbert ${\mathcal
A}$-modules, and $\varphi: {\mathcal X}\times {\mathcal Y}\to
[0,\infty)$ is a function. We start our work with the following
definition.

\begin{definition}
A (not necessarily linear) mapping $f:{\mathcal X}\to {\mathcal Y}$
is called a $\varphi$-perturbation of an adjointable mapping if
there exists a (not necessarily linear) corresponding mapping
$g:{\mathcal Y}\to {\mathcal X}$ such that
\begin{eqnarray}\label{fg}
  \|\langle f(x),y\rangle-\langle x,g(y)\rangle\|\leq \varphi(x, y)
  \qquad \big(x\in {\mathcal D}(f) \subseteq {\mathcal X}, y\in {\mathcal D}(g) \subseteq {\mathcal Y}\big)\,.
\end{eqnarray}
\end{definition}

To prove our main result, we need the following known lemma
(cf.~\cite[p.~8]{LAN}) that we prove it for the sake of
completeness.

\begin{lemma}\label{lemma}
Every densely defined adjointable mapping between Hilbert
$C^*$-modules over a fixed $C^*$-algebra ${\mathcal A}$ is
${\mathcal A}$-linear. If the adjointable mapping is everywhere
defined then it is bounded.
\end{lemma}
\begin{proof}
Let $f: {\mathcal X} \to {\mathcal Y}$ and $g: {\mathcal Y} \to
{\mathcal X}$ be a pair of densely defined adjointable mappings
between two Hilbert $C^*$-modules ${\mathcal X}$ and ${\mathcal Y}$.
For every $x_1,x_2, x_3 \in {\mathcal D}(f) \subseteq {\mathcal X}$,
every $y\in {\mathcal D}(g) \subseteq {\mathcal Y}$, every $\lambda
\in {\mathbb C}$, every $a\in {\mathcal A}$ the following equality
holds:

\begin{eqnarray*}
\langle f(\lambda x_1+x_2+x_3a),y\rangle
   & = &  \langle \lambda x_1+x_2+x_3a, g(y)\rangle \\
   & = &  \lambda \langle x_1, g(y)\rangle + \langle x_2, g(y)\rangle + a^*\langle x_3, g(y)\rangle\\
   & = &  \lambda \langle f(x_1),y\rangle + \langle f(x_2),y\rangle +a^*\langle f(x_3),y\rangle \\
   & = &  \langle \lambda f(x_1)+f(x_2)+f(x_3)a,y\rangle \, .
\end{eqnarray*}
By the density of the domain of $g$ in ${\mathcal Y}$ the equality
yields the ${\mathcal A}$-linearity of $f$.

Now, suppose $f$ and $g$ to be everywhere defined on ${\mathcal X}$
and ${\mathcal Y}$, respectively. For each $x$ in the unit sphere of
${\mathcal X}$, define $\tau_x:{\mathcal Y}\to {\mathcal A}$ by
$\tau_x(y)=\langle f(x),y\rangle=\langle x,g(y)\rangle$. Then
$\|\tau_x(y)\|\le\|x\| \|g(y)\|\le\|g(y)\|$ for any $x$ from the
unit ball. By the Banach--Steinhaus theorem we conclude that the set
$\{\|\tau_x\|: x\in {\mathcal X}, \|x\|\leq 1\}$ is bounded. ´ Due
to the equality $\|f(x)\|=\sup_{\|y\| \leq 1}\|\langle
f(x),y\rangle\|=\sup_{\|y\|=1}\|\tau_x(y)\|= \|\tau_x\|$ the mapping
$f$ has to be bounded.
\end{proof}

\begin{theorem} \label{main}
Let $f:{\mathcal X}\to {\mathcal Y}$ be a $\varphi$-perturbation of
an adjointable mapping with corresponding  mapping $g:{\mathcal
Y}\to {\mathcal X}$. Suppose that the mappings $f$ and $g$ are
everywhere defined on the respective Hilbert $C^*$-modules.
Furthermore, suppose that for some sequence $\{c_n\}$ of non-zero
complex numbers either both the conditions (\ref{phi1x}) and
(\ref{phi1y}) or both the conditions (\ref{phi2x}) and (\ref{phi2y})
below hold for the perturbation bound mapping $\varphi(x,y)$:
\begin{eqnarray}\label{phi1x}
\lim_{n\to\infty} |c_n|^{-1} \varphi(c_n x,y) = 0 \qquad(x \in
{\mathcal X}, y\in {\mathcal Y})
\end{eqnarray}
\begin{eqnarray}\label{phi1y}
\lim_{n\to\infty} |c_n|^{-1} \varphi(x,c_n y) = 0 \qquad(x \in
{\mathcal X}, y\in {\mathcal Y}) \, ,
\end{eqnarray}
\begin{eqnarray}\label{phi2x}
\lim_{n\to\infty} |c_n| \varphi(c_n^{-1}x,y) = 0 \qquad(x \in
{\mathcal X}, y\in {\mathcal Y})
\end{eqnarray}
\begin{eqnarray}\label{phi2y}
\lim_{n\to\infty} |c_n| \varphi(x,c_n^{-1}y) = 0 \qquad(x \in
{\mathcal X}, y\in {\mathcal Y}) \, .
\end{eqnarray}

Then $f$ is adjointable. In particular, $f$ is bounded, continuous
and ${\mathcal A}$-linear, as well as its adjoint is $g$.
\end{theorem}

\begin{proof}
Let $\lambda \in {\mathbb C}$ be an arbitrary number. Replacing $x$
by $\lambda x$ in (\ref{fg}), we get
\begin{eqnarray*}
\|\langle f(\lambda x),y\rangle-\langle \lambda x, g(y)\rangle\|\leq
\varphi(\lambda x,y)\, ,
\end{eqnarray*}
and since a multiplication of (\ref{fg}) by $| \lambda |$ yields
\begin{eqnarray*}
\|\langle \lambda f(x),y\rangle-\langle \lambda x, g(y)\rangle\|\leq
|\lambda|\varphi(x,y)
\end{eqnarray*}
we obtain
\begin{eqnarray}\label{lf}
\|\langle f(\lambda x),y\rangle-\langle \lambda f(x), y\rangle\|\leq
\varphi(\lambda x,y) + |\lambda|\varphi(x,y)
\end{eqnarray}
If (\ref{phi1y}) holds, we take $c_n y$ instead $y$ in (\ref{lf}) to
get
\begin{eqnarray*}
\|\langle f(\lambda x),y\rangle-\langle \lambda f(x), y\rangle\|\leq
|c_n|^{-1} \varphi(\lambda x,c_n y) + |\lambda||c_n|^{-1}
\varphi(x,c_n y)
\end{eqnarray*}
and, as $n\to \infty$, we obtain
\begin{eqnarray}\label{hom}
\langle f(\lambda x),y\rangle = \langle \lambda f(x), y\rangle
\qquad(x \in {\mathcal X}, y\in {\mathcal Y})\,.
\end{eqnarray}
If (\ref{phi2y}) holds, we take $c_n^{-1} y$ instead $y$ in
(\ref{lf}) and we arrive also at (\ref{hom}). Therefore,
\begin{eqnarray}\label{hom1}
f(\lambda x) = \lambda f(x) \qquad (x \in {\mathcal X}, \lambda \in
{\mathbb C})\,.
\end{eqnarray}
If (\ref{phi1x}) holds, we take $c_n x$ instead $x$ in (\ref{fg}) to
get
\begin{eqnarray*}
\|\langle f(c_n x),y\rangle-\langle c_n x, g(y)\rangle\|\leq
\varphi(c_n x,y)
\end{eqnarray*}
and, by (\ref{hom}), we obtain
\begin{eqnarray*}
\|\langle f(x),y\rangle-\langle x, g(y)\rangle\|\leq |c_n|^{-1}
\varphi(c_n x,y)
\end{eqnarray*}
Taking the limit as $n \to \infty$ we conclude that
\begin{eqnarray}\label{hom2}
\langle f(x),y\rangle = \langle x, g(y)\rangle \qquad(x \in
{\mathcal X}, y\in {\mathcal Y})\,.
\end{eqnarray}
Hence $f$ is adjointable and admits the mapping $g$ as its
adjoint.\\ Alternatively, if (\ref{phi2x}) holds, we take $c_n^{-1}
x$ instead $x$ in (\ref{lf}) and arrive at the same conclusion
(\ref{hom2}). By Lemma \ref{lemma} the mapping $f$ is ${\mathcal
A}$-linear and bounded with the adjoint $g$.
\end{proof}

Using the sequence $c_n=2^n$ we get the following results.

\begin{corollary}
If $f:{\mathcal X}\to {\mathcal Y}$ is an everywhere defined
$\varphi$-perturbation of an adjointable mapping, where
$\varphi(x,y)=\varepsilon\,\|x\|^p\,\|y\|^q\,\,(\alpha>0, p\neq 1,
q\neq 1)$, then $f$ is adjointable and hence a bounded $C^*$-linear
mapping.
\end{corollary}

\begin{corollary}
If $f:{\mathcal X}\to {\mathcal Y}$ is an everywhere defined
$\varphi$-perturbation of an adjointable mapping, where
$\varphi(x,y)=\varepsilon_1\,\|x\|^p\ +
\varepsilon_2\,\|y\|^q\,\,(\varepsilon_1\geq 0, \varepsilon_2 \geq
0, p\neq 1, q\neq 1)$, then $f$ is adjointable and hence a bounded
$C^*$-linear mapping.
\end{corollary}

We would like to point out that the proof of Theorem \ref{main}
works equally well in the case that the functions $f$ and $g$ are
well-defined merely on norm-dense subsets of ${\mathcal X}$ and
${\mathcal Y}$, respectively. This case covers the situation of
pairs of adjoint to each other, densely defined ${\mathcal
A}$-linear operators between pairs of Hilbert ${\mathcal
A}$-modules. However, since boundedness cannot be demonstrated, in
general, in that case we arrive at the following statement:

\begin{theorem} \label{main2}
Let $f:{\mathcal X}\to {\mathcal Y}$ be a $\varphi$-perturbation of
an adjointable mapping with corresponding mapping $g:{\mathcal Y}\to
{\mathcal X}$. Suppose, that the mappings $f$ and $g$ are densely
defined on the respective Hilbert $C^*$-modules. Furthermore,
suppose that for the perturbation bound mapping $\varphi(x,y)$
either both the conditions (\ref{phi1x}) and (\ref{phi1y}), or both
the conditions (\ref{phi2x}) and (\ref{phi2y}) hold. Then $f$ is
adjointable. In particular, $f$ is ${\mathcal A}$-linear, as well as
its adjoint is $g$.
\end{theorem}

\begin{corollary}
The equation $f(x)^*y=xg(y)^* \quad (x\in {\mathcal I}, y\in
{\mathcal J})$ is superstable, where \newline $f:{\mathcal I}\to
{\mathcal J}$ and $g:{\mathcal J}\to {\mathcal I}$ are adjoint to
each other, densely defined ${\mathcal A}$-linear mappings between
right ideals ${\mathcal I}, {\mathcal J}$ of ${\mathcal A}$.
\end{corollary}

The critical case of $\varphi$-perturbations is that one where the
function $\varphi$ satisfies neither the pair of conditions (i) and
(ii), nor the pair of conditions (i') and (ii'). We demonstrate that
there may exist $\varphi$-perturbed bounded $C^*$-linear mappings
$f$ on certain types of Hilbert $C^*$-modules ${\mathcal X}$ over
certain $C^*$-algebras ${\mathcal A}$ which are not adjointable.
Moreover, any non-adjointable bounded $C^*$-linear mapping $f$ on
suitably selected Hilbert $C^*$-modules ${\mathcal X}$ can be
$\varphi$-perturbed by ``compact'' operators on ${\mathcal X}$ using
this type of perturbation functions.

\begin{proposition}
Let ${\mathcal X}$ be a Hilbert ${\mathcal A}$-module over a given
$C^*$-algebra ${\mathcal A}$. Suppose there exists a non-adjointable
bounded ${\mathcal A}$-linear mapping $f: {\mathcal X} \to {\mathcal
X}$, (so ${\mathcal X}$ cannot be self-dual by \cite{LAN, FRA}).
Then there exist (at least countably many) positive constants
$c_\alpha$ and respective ``compact'' ${\mathcal A}$-linear
operators $K_\alpha: {\mathcal X} \to {\mathcal X}$ ($\alpha \in I$)
such that $f$ is $\phi$-perturbed for a function $\phi(x,y)=
c_\alpha \cdot \|x\| \cdot \|y\|$ and for $g=K_\alpha^*$.
\end{proposition}

\begin{proof}
By results of Huaxin~Lin \cite{LIN1}, \cite[Thm. 1.5]{LIN2}, the
Banach algebra $End_{\mathcal A}({\mathcal X})$ of all bounded
${\mathcal A}$-linear mappings on ${\mathcal X}$ is the left
multiplier algebra of the C*-algebra $K_{\mathcal A}({\mathcal X})$
of all ``compact'' ${\mathcal A}$-linear operators on ${\mathcal
X}$. Since $End_{\mathcal A}({\mathcal X})$ is the completion of
$K_{\mathcal A}({\mathcal X})$ with respect to the left strict
topology defined by the set of semi-norms $\{ \|\cdot K \| : K \in
K_{\mathcal A} ({\mathcal X}) \}$, there exists a bounded net $\{
K_\alpha : \alpha \in I \}$ of ``compact'' operators such that the
set $\{ K_\alpha K : \alpha \in I \}$ converges with respect to the
operator norm to $fK$ for any single ``compact'' operator $K$.
Therefore,
\[
0= \lim_{\alpha \in I} \| \langle (fK-K_\alpha K)(x),y \rangle \|
  = \lim_{\alpha \in I} \| \langle (f-K_\alpha)K(x),y \rangle \|
\]
for any ``compact'' operator $K$. However, the set $\{ K(x) : K \in
K_{\mathcal A}({\mathcal X}), x \in {\mathcal X} \}$ is norm-dense
in ${\mathcal X}$, hence
\[
       \| \langle f(x),y \rangle - \langle K_\alpha(x),y \rangle \|
\leq \|f-K_\alpha\| \cdot \|x\| \cdot \|y\|
\]
for any $x,y \in {\mathcal X}$ and any $\alpha \in I$. Setting
$c_\alpha = \|f-K_\alpha\|$ for any fixed index $\alpha$ and taking
into account the adjointability of the operators $\{ K_\alpha \}$ we
arrive at the desired result.
\end{proof}

\begin{corollary}
Let ${\mathcal X}$ be a Hilbert ${\mathcal A}$-module over a given
$C^*$-algebra ${\mathcal A}$. Suppose there exists a non-adjointable
bounded ${\mathcal A}$-linear mapping $f: {\mathcal X} \to {\mathcal
X}$. Then there does not exist any $\varphi$-perturbation of $f$
such that $\varphi(x,y)$ satisfies either both the conditions
(\ref{phi1x}) and (\ref{phi1y}) or both the conditions (\ref{phi2x})
and (\ref{phi2y}).
\end{corollary}

\vspace{1cc}

{\small

\noindent Michael Frank: Hochschule f\"ur Technik, Wirtschaft und
Kultur (HTWK) Leipzig, Fachbereich Informatik, Mathematik und
Naturwissenschaften (FbIMN), PF 301166, D-04251 Leipzig, Germany.\\
mfrank@imn.htwk-leipzig.de \vspace{0.8cc}

\noindent Pasc G\u avru\c ta: Department of Mathematics, University
`Politehnica' of Timi\c{s}oara, Piata Victoriei, No. 2, 300006
Timi\c{s}oara, Romania.\\
pgavruta@yahoo.com \vspace{0.8cc}

\noindent Mohammad Sal Moslehian: Department of Pure Mathematics,
Ferdowsi  University of Mashhad, P. O. Box 1159, Mashhad 91775, Iran;\\
Centre of Excellence in Analysis on Algebraic Structures (CEAAS),
Ferdowsi University of Mashhad, Iran.\\
moslehian@ferdowsi.um.ac.ir and moslehian@ams.org \vspace{0.8cc}
}
\end{document}